\pdfoutput=1
\RequirePackage{silence}
\documentclass[a4paper,svgnames]{amsart}
\usepackage[hmarginratio={1:1},vmarginratio={1:1},lmargin=61.0pt,tmargin=60.0pt]{geometry}
\usepackage{float}

\usepackage{booktabs}
\allowdisplaybreaks

\usepackage[T1]{fontenc}
\usepackage{latexsym,exscale,amsfonts,amssymb,mathtools}
\usepackage{amsmath,amsthm,amsfonts,amssymb,amscd,textcomp,bbm}
\usepackage{mathrsfs,stackrel}
\usepackage{etoolbox}
\usepackage{tikz,tikz-cd}
\usetikzlibrary{matrix,arrows.meta,decorations.markings,shapes.multipart,decorations.pathreplacing,decorations.pathmorphing,shapes.geometric}
\usepackage{aliascnt}
\usepackage{ytableau}
\usepackage{xparse}
\usepackage{cite}

\usepackage{dynkin-diagrams}
\usepackage{nicematrix}
\usepackage{tikz-3dplot}
\usepackage{mathrsfs}
\DeclareMathAlphabet{\mathscrbf}{OMS}{mdugm}{b}{n}
\DeclareMathOperator{\M}{M}
\DeclareMathOperator{\GL}{GL}

\usepackage{array}
\newcolumntype{C}{>{$}c<{$}}
\newcolumntype{P}[1]{>{\centering\arraybackslash}p{#1}} 
\DeclarePairedDelimiter\floor{\lfloor}{\rfloor}

\usepackage{xcolor,colortbl}

\definecolor{mygray}{gray}{0.6}
\definecolor{mygraydark}{gray}{0.4}
\definecolor{mygraylight}{gray}{0.85}
\definecolor{spinach}{RGB}{46,139,87}
\definecolor{tomato}{RGB}{255,99,71}
\definecolor{orchid}{RGB}{143,40,194}
\definecolor{neon}{RGB}{77,77,255}
\definecolor{lightneon}{RGB}{110,110,255}
\definecolor{pumpkin}{RGB}{224,180,80}
\definecolor{citron}{RGB}{190,180,90}

\definecolor{lava}{RGB}{207,16,32}
\definecolor{cream}{RGB}{255,253,208}
\definecolor{verdigris}{RGB}{67,179,174}
\definecolor{Black}{RGB}{0,0,0}
\definecolor{mydarkblue}{RGB}{10,10,170}
\definecolor{darkspinach}{RGB}{20,70,20}
\definecolor{darktomato}{RGB}{155,40,30}
\definecolor{darkorchid}{RGB}{50,10,100}
\definecolor{darklava}{RGB}{150,8,16}

\definecolor{zero}{RGB}{0,0,0}
\definecolor{one}{RGB}{255,0,0}
\definecolor{two}{RGB}{0,255,0}
\definecolor{three}{RGB}{0,0,255}

\usepackage{enumitem}
\setlist[enumerate]{itemsep=0.15cm,label=\emph{\upshape\arabic*.}}
\setlist[enumerate,2]{itemsep=0.15cm,label=\emph{\upshape(\roman*)}}
\setlist[enumerate,3]{itemsep=0.15cm,label=\emph{\upshape(\Alph*)}}

\newcommand{\Z}{\mathbb{Z}}
\newcommand{\Q}{\mathbb{Q}}

\newcommand{\C}{\mathbb{C}}

\DeclareMathOperator{\Ind}{Ind}
\DeclareMathOperator{\Res}{Res}

\usepackage{chngcntr}
\counterwithout{equation}{section}     

\def\equationautorefname~#1\null{(#1)\null}

\theoremstyle{plain}
\newtheorem{Proposition}{Proposition}
\newtheorem{Theorem}{Theorem}

\newtheorem{Lemma}{Lemma}
\newtheorem{Conjecture}{Conjecture}

\theoremstyle{definition}

\newtheorem{Example}{Example}

\theoremstyle{remark}
\newtheorem{Remark}{Remark}

\AtEndEnvironment{Corollary}{\null\hfill$\square$}%
\AtEndEnvironment{Definition}{\null\hfill$\Diamond$}%
\AtEndEnvironment{ClassificationProblem}{\null\hfill$\Diamond$}%
\AtEndEnvironment{Notation}{\null\hfill$\Diamond$}%
\AtEndEnvironment{Example}{\null\hfill$\Diamond$}%
\AtEndEnvironment{Examples}{\vskip-10mm\null\hfill$\Diamond$}%
\AtEndEnvironment{Remark}{\null\hfill$\Diamond$}%
\AtEndEnvironment{Question}{\null\hfill$\Diamond$}%

\usepackage[hypertexnames=false]{hyperref}
\usepackage{bookmark}
\hypersetup{
pdftoolbar=true,
pdfmenubar=true,
pdffitwindow=false,
pdfstartview={FitH},
pdftitle={Growth problems for representations of finite monoids},
pdfauthor={David He and Daniel Tubbenhauer},
pdfsubject={},
pdfcreator={David He and Daniel Tubbenhauer},
pdfproducer={David He and Daniel Tubbenhauer},
pdfkeywords={},
pdfnewwindow=true,
colorlinks=true,
linkcolor=mydarkblue,
citecolor=teal,
filecolor=magenta,
urlcolor=orchid,
linkbordercolor=lava,
citebordercolor=teal,
urlbordercolor=orchid,
linktocpage=true
}

\def\makeautorefname#1#2{\csdef{#1autorefname}{#2}}
\makeautorefname{section}{Section}%
\makeautorefname{subsection}{Section}%
\makeautorefname{subsubsection}{Section}%

\begin{document}

\title{Growth problems for representations of finite monoids}
\author{David He and Daniel Tubbenhauer}
\address{D. H.: The University of Sydney, School of Mathematics and Statistics, Australia}
\email{d.he@sydney.edu.au}

\address{D.T.: The University of Sydney, School of Mathematics and Statistics F07, Office Carslaw 827, NSW 2006, Australia, \href{http://www.dtubbenhauer.com}{www.dtubbenhauer.com}, \href{https://orcid.org/0000-0001-7265-5047}{ORCID 0000-0001-7265-5047}}
\email{daniel.tubbenhauer@sydney.edu.au}

\begin{abstract}
We give a conjecture for the asymptotic growth rate of the number of indecomposable summands in the tensor powers of representations of finite monoids, expressing it in terms of the (Brauer) character table of the monoid's group of units. 
We prove it under an additional hypothesis. We also give (exact and asymptotic) formulas for the growth rate of the length of the tensor powers when working over a good characteristic. As examples, we compute the growth rates for the full transformation monoid, the symmetric inverse monoid, and the monoid of 2 by 2 matrices. We also provide code used for our calculation.
\end{abstract}


\addtocontents{toc}{\protect\setcounter{tocdepth}{1}}

\maketitle

\tableofcontents

\section{Introduction}
\subsection{Growth problems}
Throughout let $M$ be a finite monoid, and let $k$ be a splitting field for $M$ of characteristic $p\ge 0$. For a $kM$-module $V$, we are interested in the quantity
$$b(n)=b^{M,V}(n)=\text{\# $M$-indecomposable summands in} \ V^{\otimes n } \ \text{(counted with multiplicity)}.$$
In particular, we ask the following questions:
\begin{enumerate}[label={(\arabic*)}]
    \item \label{1} Can we find a formula for $b(n)$, or more generally a formula for an asymptotic expression $a(n)$ with $b(n)\sim a(n)?$ (Here $b(n)\sim a(n)$ if they are asymptotically equal: ${b(n)}/{a(n)}\xrightarrow{ n \to \infty } 1.$)
    \item \label{2} Can we understand the \textit{rate of geometric convergence} by quantifying how fast $\lvert b(n)/a(n)-1\rvert$ converges to 0?
    \item \label{3} Similarly, can we bound the \textit{variance} $\lvert b(n)-a(n)\rvert$?
\end{enumerate}
The questions \ref{1}--\ref{3} above are a special case of \textit{growth problems}, which following \cite{lacabanne2024asymptotics} we may define for any additive Krull--Schmidt monoidal category. The growth problems and related questions for various categories have been recently studied in the papers \cite{coulembier2023asymptotic, coulembier2024fractalbehaviortensorpowers,he2024growthproblemsrepresentationsfinite,KhSiTu-monoidal-cryptography,lacabanne2023asymptotics, lacabanne2024asymptotics,lachowska2024tensorpowersvectorrepresentation, larsen2024boundsmathrmsl2indecomposablestensorpowers}. In particular, if our monoid $M$ is a group, then the situation is well-understood from \cite[Theorem 1]{he2024growthproblemsrepresentationsfinite}, which gives an expression for the asymptotic formula $a(n)$ in terms of the (Brauer) character table. A main goal of this paper is to generalize this result to an arbitrary finite monoid. 

In addition to studying $b(n)$, which counts the number of indecomposable summands of $V^{\otimes n}$, we also ask the analogue of the questions \ref{1}--\ref{3} for the related quantity $l(n)=l(V^{\otimes n})$, where $l(W)$ denotes the \textit{length} of $W$, \textit{i.e.} the number of composition factors. If $kM$ is semisimple, then $b(n)=l(n)$.  

\subsection{A conjecture and evidence}\label{conj section}

Let $G$ denote the group of units of $M$. Let $Z_V(G)$ denote the set of elements $g$ in $G$ which act as scalars on $V$, and denote the corresponding scalars by $\omega_V(g)$. Let $g_t, 1\le t \le N$ be a complete set of representatives for the $p$-regular conjugacy classes of $G$ (if $p=0$, these are all the conjugacy classes).  If $V$ is a $kM$-module, recall from \cite[Section 2]{lacabanne2024asymptotics} that the corresponding \textit{fusion graph} $\Gamma$ is the (oriented and weighted) graph whose vertices are indecomposable $kM$-modules which are summands of $V^{\otimes n}$, for some $n\ge 0$, and that there is an edge of weight $m$ from the vertex $V_i$ to the vertex $V_j$ if $V_i$ occurs $m$ times in the direct sum decomposition of $V\otimes V_j$. The corresponding (potentially countably infinite) adjacency matrix is called the \textit{action matrix}. If $V$ is a $kM$-module, let $\Res_G(V)$ denote the $kG$-module that comes from restricting the action on $V$. Let $\lambda^{\mathrm{sec}}$ denote any second largest eigenvalue (in terms of modulus) of the action matrix of $V$. We use the usual capital $O$ notation. We conjecture that the asymptotic growth rate $a(n)$ is the same as that of $\Res_G(V)$. \textcolor{black}{We note that if $M=G$ is a group, \autoref{conj} recovers the statement of \cite[Theorem 1]{he2024growthproblemsrepresentationsfinite}.}
\begin{Conjecture}\label{conj}
Suppose $V$ is such that no element apart from $1\in M$ acts as identity, then

\begin{enumerate}
    \item The asymptotic formula is \begin{equation}\label{eqn:asym} a(n)=\frac{1}{\lvert G\rvert}\sum_{\substack{1\le t\le N \\ g_t \in Z_V(G)}} S_{t}\big(\omega_V(g_t)\big)^n\cdot (\dim V)^n, 
    \end{equation}
    where $S_t$ is the sum over entries of the column corresponding to $g_t^{-1}$ in the (irreducible) Brauer character table. In other words, $V$ has the same asymptotic growth rate as the $kG$-module $\Res_G(V)$, cf. \cite[Theorem 1]{he2024growthproblemsrepresentationsfinite}.
\item $\lvert b(n)/a(n)-1\rvert \in \mathcal{O}(\lvert \lambda^{\mathrm{sec}}/\dim V \rvert^n + n^{-c})$, for some constant $c>0$, and
    \item $\lvert b(n)-a(n) \rvert \in \mathcal{O}(\lvert \lambda^{\mathrm{sec}}\rvert^n + n^{d})$ for some constant $d > 0$.
\end{enumerate}
\end{Conjecture}
\begin{Remark}
When working with group representations $V$, we generally require that $V$ be faithful. The requirement that no elements apart from $1\in M$ acts as identity should be seen as the replacement of this requirement in the more general setting.     
\end{Remark}
\textcolor{black}{The reason we expect \autoref{conj} to be true is the following: The collection $\mathcal{G}$ of indecomposable $kM$-modules which are $kG$-modules (that is, $M\setminus G$ acts as 0) forms a \textit{tensor ideal} in the sense that if $K \in \mathcal{G}$ and $V$ is any $kM$-module, then $V\otimes K$ is a direct sum of modules in $\mathcal{G}$. By the technical machinery introduced in \cite{lacabanne2024asymptotics}, our conjecture that $\mathcal{G}$ determine the asymptotic growth rate is roughly equivalent to the conjecture that every vertex in the fusion graph has a path to $\mathcal{G}$, which we suspect is true based on experimental evidence.}

As evidence for \autoref{conj}, over fixed finite fields $k$ (\textit{e.g. $\mathbb{F}_{11}$}) we computed the growth rates of all projective indecomposable $kM$-modules, where $M$ is either a monoid of order $\le 7$ (34129 of these) or is the direct product $T_3\times L$ where $T_3$ is the full transformation monoid on 3 elements, and $L$ is a monoid of order $\le 5$ (273 of these). We were not able to find a counterexample to the formula in \autoref{eqn:asym}. To compute the growth rates, we used GAP's 
\verb|Smallsemi| package to obtain the multiplication table for the monoids of small order, and from these we constructed the monoid algebras in Magma as matrix algebras. Finally, we computed the action matrix as in the appendix of \cite{he2024growthproblemsrepresentationsfinite}. We performed the calculations over finite fields as Magma's Meataxe algorithm for decomposing modules is much slower over $\Q$.

\begin{Remark}
The code used in this paper is available on the GitHub repository \cite{code}.  
\end{Remark}

\subsection{Main results}
In \autoref{Main}, we prove \autoref{conj} under the hypothesis that some projective $kG$-module is injective as a $kM$-module. While the hypothesis is somewhat restrictive, it is satisfied by important classes of monoids such as the full transformation monoids $T_m$ (in characteristic 0), the monoid of $2\times 2$ matrices $\M(2,q)$ (over defining characteristic), and all monoids with semisimple monoid algebra. We compute asymptotic formulas for these examples in \autoref{examples}. In addition, we obtain general (exact and asymptotic) formulas for $l(n)$ when the characteristic of $k$ does not divide the order of any of $M$'s maximal subgroups. As $b(n)=l(n)$ when $kM$ is semisimple, we obtain a formula for $b(n)$ in this case, which generalizes \cite[Theorem 3]{he2024growthproblemsrepresentationsfinite}.

\noindent\textbf{Acknowledgments.} We would like to express our gratitude to Volodymyr (Walter) Mazorchuk for insightful email exchanges on the topic of monoid representations. In particular, the observation discussed in \autoref{SS:Counterexample} is credited to Walter. \textcolor{black}{We would also like to thank the referee for giving the paper a careful reading and providing many helpful comments.}

DT acknowledges the support of the ARC Future Fellowship (FT230100489).
Lastly, DT has reached a state of complete insanity.

\section{Growth rate in general}

We say a function $f$ satisfies $f\in\Theta^{\prime}(g)$ if there exists a constant $A\in\mathbb{R}_{>0}$ such that $A\cdot g(n)\leq f(n)\leq g(n)$ for 
all $n>n_{0}$ for some fixed $n_{0}\in\mathbb{N}$. We can say the following about a general semigroup representation:

\begin{Proposition}\label{p:finitegroup}
Let $S$ be a finite semigroup, and $V$ be a finite dimensional $kS$-module.
Then
\begin{gather*}
b(n),l(n)\in\Theta^{\prime}\big((\dim V)^{n}\big).
\end{gather*}
\end{Proposition}

\begin{proof}
The upper bound is immediate, and it remains to justify the lower bound. Moreover, $b(n)\leq l(n)$, so it remains to justify the lower bound for $b(n)$.

To this end, let $\M(m,k)$ denote the monoid of $m\times m$ matrices with values in $k$.
Let $\GL(m,k)$ denote the invertible matrices in $\M(m,k)$. As in \cite{Coulembier_2023}, $\mathcal{O}(M(m,k))$ is a subcoalgebra in $\mathcal{O}(\GL(m,k))$ and the category of $k\M(m,k)$-modules can be identified with the category of polynomial 
$k\GL(m,k)$-modules, so the number of summands over $\M(m,k)$ is the same as the number of summands over $\GL(m,k)$. 

Now, any finite semigroup $S$ with $\#S=m$ can be realized in $\M(m,\mathbb{F}_{2})$ 
(via the regular representation over $\mathbb{F}_{2}$), so we get a bound of $b(n)$ from below by $b(n)$ for $\M(m,\mathbb{F}_{2})$ which is the same as for 
$\GL(m,\mathbb{F}_{2})$, by the above. 
The desired lower bound now follows from \cite[Proposition 2.2]{Coulembier_2023}.
\end{proof}

We focus solely on monoids for the remainder of the paper. However, the following example demonstrates that for general semigroups, there might be a final basic class represented by the null representation (which does not exist for a monoid), allowing the application of standard theory, as discussed in \cite{lacabanne2024asymptotics}.

\begin{Example}
Let $S$ be the semigroup with underlying set $G_1\cup G_2\cup \{0\}$, with multiplication such that $G_1\cong \Z/2\Z$, $G_2\cong S_3$ and $G_1G_2=\{0\}$. The monoid algebra $\C S$ is semisimple with 7 simple modules: the trivial module, the one-dimensional null representation $Z$ on which $S$ acts as 0, and the simple modules for $G_1\cong \Z/2\Z$ and $G_2\cong S_3$ respectively, extended so that elements outside $G_1$ (resp. $G_2$) in $S$ act as 0. Let $X$ be the $\C S$-module corresponding to the non-trivial one-dimensional simple $\C G_1$-module, and let $Y$ be the $\C S$-module corresponding to the two-dimensional simple `standard representation' of $\C G_2$. If we set $V=Y\oplus X$, then the vertex $Z$ is a sink node in the fusion graph of $V$ which determines the asymptotic growth rate. We note that whenever $Z$ shows up in the graph, it is a final basic class (FBC) in the sense of \cite[Definition 5.5]{lacabanne2024asymptotics} and we must have $a(n)=(\dim V)^n$. 

\begin{figure}[ht]
    \centering
    \begin{minipage}{0.40\textwidth}
    \centering
    \includegraphics[width=0.6\textwidth]{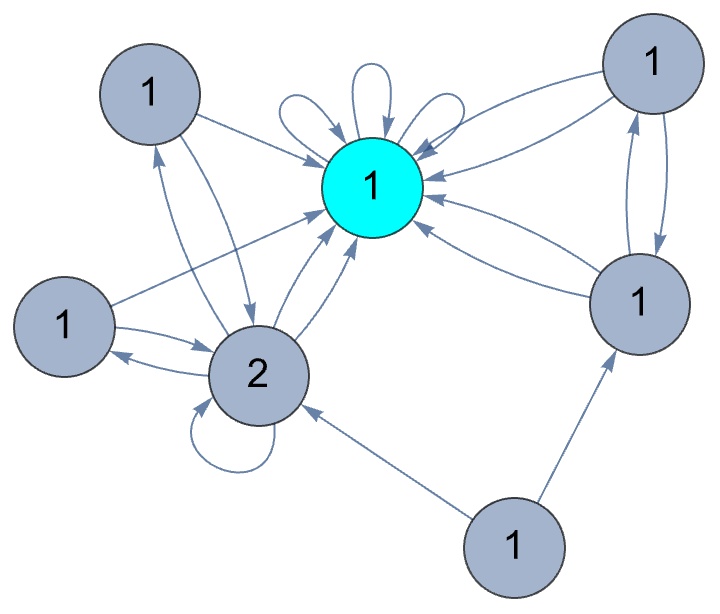} 
    \end{minipage}
    \begin{minipage}{0.4\textwidth}
     \centering
     
\[
\begin{pmatrix}
0 & 0 & 0 & 0 & 0 & 0 & 0 \\
1 & 0 & 0 & 0 & 1 & 0 & 0 \\
1 & 0 & 1 & 0 & 0 & 1 & 1 \\
0 & 2 & 2 & 3 & 2 & 1 & 1 \\
0 & 1 & 0 & 0 & 0 & 0 & 0\\
0 & 0 & 1 & 0 &0 & 0 & 0 \\
0 & 0 & 1 & 0 & 0&0 &0 
\end{pmatrix}
\]

    \end{minipage}
    \caption{The fusion graph and action matrix for $V=Y\oplus X$. The vertex colored in cyan is the null representation. In this case $a(n)=3^n$.} 
    \label{fig:semigroupsumexample}
\end{figure}

\end{Example}

\section{Main results}

\subsection{Counting summands} \label{main section}

 \textcolor{black}{This section focuses on the study of $b(n)$.} As before, if $V$ is a $kM$-module let $\Res_G(V)$ denote the $kG$-module that comes from restriction of $V$. If $W$ is a $kG$-module, let $\Ind_G(W)$ denote the induced $kM$-module on which elements in $M\setminus G$ act as 0. We write $P(W)$ for the projective cover of $W$. We note that induction and restriction give inverse equivalences between the category of $kG$-modules and the full subcategory of $kM$-modules with all nonunits acting as 0. Let $\Gamma^G$ denote the subgraph in the fusion graph $\Gamma$ for $V$ induced by the vertices $\Ind_G(P(W))$ where $W$ is an irreducible $kG$-module.

Let $\Gamma_G^P$ denote the \textit{projective cell} in the fusion graph for $\Res_G(V)$, then $\Gamma_G^P$ and $\Gamma^G$ are isomorphic as graphs, as we have 
\begin{equation} \label{tensor}
V\otimes \Ind_G(P(W))= \Ind_G(\Res_G(V)\otimes P(W)).
\end{equation} \textcolor{black}{Recall from \cite[\S 4]{lacabanne2024asymptotics} that a \textit{final basic class} (FBC) of a fusion graph is a strongly-connected component whose PF dimension (= spectral radius of the adjacency submatrix in the case when it is finite, the only case we need) is maximal among all strongly-connected components, and moreover has no path to any other strongly-connected component with  maximal PF dimension. If there is no path leaving the FBC at all, we say it is final in the whole graph.} The following theorem establishes \autoref{conj} under the assumption that some module in $\Gamma^G$ is injective. 

\begin{Theorem} \label{Main}
Suppose that $V$ is a $kM$-module such that \begin{enumerate}[label=\roman*)]
    \item No element apart from $1 \in M$ acts as the identity, and
    \item There is an irreducible $kG$-module $W$ such that $\Ind_G(P(W))$ is an injective $kM$-module,
\end{enumerate}
then the statements of \autoref{conj} hold.

\end{Theorem}
\begin{proof}
Let $\Gamma^G$ be the subgraph of $\Gamma$ as above, then it is strongly connected since $\Res_G(V)$ is faithful by assumption i), and the \textit{projective cell} in the fusion graph of a faithful $kG$-module is strongly-connected (cf. the proof of \cite[Proposition 4.22]{lacabanne2024asymptotics}). \textcolor{black}{We will first show that the growth problem for $V$ is \textit{sustainably positive recurrent} in the sense of \cite[Definition 5.5]{lacabanne2024asymptotics}, with $\Gamma^G$ being the unique FBC. This allows us to invoke \cite[Theorem 5.10]{lacabanne2024asymptotics} from which the statements of the theorem easily follow}. 

We can see that $\Gamma^G$ is a basic class because it has PF dimension equal to $\dim V$ (since $\Gamma^G \cong \Gamma_G^P$ and PF$\dim \Gamma_G^p=\dim \Res_G(V)=\dim V$), which is the maximal possible since $b(n)$ grows as some multiple of $(\dim V)^n$. It is moreover final (in fact, final in $\Gamma$) by \autoref{tensor} and the fact that the tensor product of any $kG$-module with a projective $kG$-module is projective.  

Next, we claim that under the assumptions of the theorem any vertex $Z$ in $\Gamma$ has a path to $\Gamma^G$. For $n \ge 0$, the set $V_n\subseteq V^{\otimes n}$ of all vectors $v\in V^{\otimes n}$ on which $M\setminus G$ acts as zero is a submodule. To see this, note that $V_n$ is clearly a vector subspace (it contains $0$, so is nonempty), and if $v\in V_n$, for all $m \in M\setminus G$ and $n\in M$ we have $m(nv)=0$ because $mn\in M\setminus G$, and so $nv\in V_n$. By assumption i), $V_n$ is nonzero for $n$ large enough, say $n=|M\setminus G|$, since for any $m \in M\setminus G$ acts non-invertibly and we can find $v_m \in V$ such that $m\cdot v_m=0$, and then the basic tensor over all $v_m$'s is in $V_n$. Thus, for some $n$ large enough, $V^{\otimes n}$ contains a nonzero submodule $V'=V_n$ on which all non-units act as zero. Because $\Res_G(V)$ is faithful, any projective indecomposable $kG$-module $P$ appears in some tensor power of $\Res_G(V)$ (see \cite[Theorem 2]{bryant1972tensor}), say $(\Res_G(V))^{\otimes k}$. Now $\Res_G(V'\otimes V^{\otimes k})=\Res_G(V')\otimes (\Res_G(V))^{\otimes k}$ has a projective indecomposable summand, and so $ V'\otimes V^{\otimes k} \subseteq V^{\otimes (k+n)}$ has a summand of the form $\Ind_G(P(W))$ for some irreducible $kG$-module $W$ (all nonunits act as zero on $V'\otimes V^{\otimes k}$, so we recover the same module by restriction followed by induction). By assumption ii), we may take $\Ind_G(P(W))$ to be injective, and so it is a summand of $V^{\otimes (k+n)}$. This shows that the trivial module $k$ has a path to $\Gamma^G$; \textcolor{black} {since $Z\cong Z\otimes k$ and $\Gamma^G$ is final in $\Gamma$, in fact any indecomposable $kM$-module $Z$ has a path to $\Gamma^G$ (we know for $N$ large enough, $k\otimes V^{\otimes N}$ has a summand in $\Gamma^G$, so $Z\otimes k \otimes V^{\otimes N}$ does too).}

Now it follows, by the same proof as \cite[Proposition 5.7]{lacabanne2024asymptotics}, that the growth problem for $V$ is sustainably positive recurrent with the FBC also final in $\Gamma$.  Statements (2) and (3) of \autoref{conj} now follow from \cite[Theorem 5.10]{lacabanne2024asymptotics}. The first statement follows by the same proof as \cite[Theorem 1]{he2024growthproblemsrepresentationsfinite}, after adapting the statement about left eigenvectors in \cite[Lemma 7]{he2024growthproblemsrepresentationsfinite} by using (Brauer) characters of monoids rather than characters of groups (the right eigenvectors are the same). For a definition of the Brauer characters of monoids we refer to \cite[\S 5]{steinberg2023modularrepresentationtheorymonoids}. 
\end{proof}

\begin{Remark} \
\begin{enumerate}
    \item When $kM$ is semisimple, condition i) in \autoref{Main} is equivalent to the requirement that no generalized conjugacy class in $M$ apart from that of the identity has character value $\dim V$. Condition ii) is automatically satisfied in this case. We will be able to say much more if $kM$ semisimple, see \autoref{semisimple} below.
    \item If $\textup{char }k \nmid \lvert G\rvert$, then $\Ind_G(P(W))=\Ind_G(W)$ for any irreducible $kG$-module $W$.
    \item We require condition ii) so that, assuming also condition i) holds, $V^{\otimes n}$ will contain a summand in $\Gamma^G$ for $n$ large enough. We suspect that this is always true, which would imply \autoref{conj}, but we do not know how to prove it in general. See also \autoref{SS:Counterexample} below.
\end{enumerate}

\end{Remark}

\subsection{Counting lengths}
\textcolor{black}{We now turn to study the related statistic of $l(n)$.} 
Let $V$ be a $kM$-module. To compute $l(n)$ we count the multiplicities now in the Grothendieck ring $G_0(kM)$ defined via short exact sequences, which coincides with the additive Grothendieck ring if $kM$ is semisimple (in this case $l(n)=b(n)$). 

Let $e_1,\dots, e_m$ denote a complete set of idempotent representatives for the regular $\mathcal{J}$-classes of $M$, and let $G_{e_1},\dots, G_{e_m}$ denote the corresponding maximal subgroups. We assume from now that $p\nmid \lvert G_{e_i}\rvert$ for all $1\le i \le m$. For a $kM$-module $V$, let $e_iV$ denote the $kG_{e_i}$-module defined via restriction. Recall from \cite[Theorem 6.5]{steinberg2016representation} that the map $\textup{Res}: G_0(kM)\to \prod_{i=1}^m G_0(kG_{e_i})$ given by $\textup{Res}([V])=([e_1V],\dots, [e_mV])$
is a ring isomorphism. Let $L$ denote the matrix corresponding to $\textup{Res}$, called the \textit{decomposition matrix}, with respect to the bases $\textup{Irr}_k(M)$ and $\bigcup_{i=1}^m \textup{Irr}_k(G_{e_i})$ respectively. Moreover, the character table $X(M)$ of $M$ is the transpose of the change of basis matrix between the basis of irreducible characters and the basis of indicator functions on the generalized conjugacy classes, and may be obtained as $X(M)=L^TX$ where $X$ is the block diagonal matrix with the blocks being the character tables of $G_{e_i}, 1\le i \le m$ (see \cite[Corollary 7.17]{steinberg2016representation}). 

If $V$ is a $kM$-module, let $N(V)$ denote the analogue of action matrix for $l(n)$. That is, with respect to the basis of isomorphism classes of irreducible $kM$-modules $\textup{Irr}_k(M)= \{S_1,\dots, S_z\}$, the $(i,j)$-th entry of $N(V)$ is the composition multiplicity $[V\otimes S_j: S_i].$ A right eigenvector of $N(e_iV)$ may be identified with an element $\mathbf{x}$ in the Grothendieck ring $G_0(kG_{e_i})$, considered as a vector with respect to the basis $\textup{Irr}_k(G_{e_i})$. By naturally identifying $G_0(kG_{e_i})$ with a subset of $\prod_{i=1}^m G_0(kG_{e_i})$, $\textbf{x}$ can be identified with a longer vector $\tilde{\textbf{x}}$ with respect to the basis $\bigcup_{i=1}^m \textup{Irr}_k(G_{e_i})$ which is supported on that subset.    

\begin{Lemma} \label{eigenvector lemma}
Let $V$ be a $kM$ module with character $\chi$ and assume $p\nmid \lvert G_{e_i}\rvert,1\le i \le m$.
\begin{enumerate}
    \item The right eigenvectors of $N(V)$ are columns of $L^{-1}\overline{X}$, where $X$ is the block diagonal matrix of character tables, with eigenvalues the character values at the corresponding generalized conjugacy classes.  
    \item The left eigenvectors of $N(V)$ are columns of the character table $X(M)$ of $M$, with eigenvalues the corresponding character values. 
\end{enumerate}
\end{Lemma}

\begin{proof}
First we prove a). For $1\le i \le m$, let $G_0(kM)_i$ denote the inverse image of $G_0(kG_{e_i})$ (identified as a subset of $\prod_{j=1}^nG_0(kG_{e_j})$) under the restriction map. That is, it is the ideal in $G_0(kM)$ consisting of elements $[W]$ such that $[e_jW]=0$ for $j\neq i$. We have the commutative diagram  \[\begin{tikzcd}
G_0(kM)_i \arrow[r,"\textup{Res}"] \arrow[d, swap, "N(V) \cdot"]  & G_0(kG_{e_i}) \arrow[d,"N(e_iV) \cdot"] \\
G_0(kM)_i \arrow[r,"\textup{Res}"] & G_0(kG_{e_i})
\end{tikzcd} \begin{tikzcd}
 \left[Y\right] \arrow[mapsto, r,"L"] \arrow[d, swap, "N(V) \cdot"] & \left[e_iY\right] \arrow[mapsto, d, "N(e_iV) \cdot"] \\
 \left[V\otimes Y\right]  \arrow[mapsto, r,"L"] & \left[e_i(V\otimes Y)\right]=\left [ e_iV\otimes e_iY\right ]
\end{tikzcd}\]   
where we think of $G_0(kG_{e_i})$ as a subset of $ \prod_{j=1}^nG_0(kG_{e_j})$ and of the elements in the Grothendieck rings as column vectors. It follows that if $\textbf{x}$ is a right eigenvector for $M(e_iV)$, $L^{-1}\tilde{\textbf{x}}$ is a right eigenvector for $N(V)$ with the same eigenvalue. By \cite[Theorem 3]{he2024growthproblemsrepresentationsfinite}, since $p\nmid \lvert G_{e_i}\rvert$ for $1\le i \le m$, $\mathbf{x}$ corresponds to the complex conjugate of a column of the character table of some $G_{e_i}$, so the claim follows.  

For b), note that $X(M)$ is by definition the transpose of the change of basis matrix between the basis of irreducible characters and the basis of indicator functions. If $\textup{Irr}_k(M)=\{S_1,\dots, S_z\}$, we have $$X(M)^T \begin{pmatrix}
  [V\otimes S_j:S_1] \\
  \vdots \\
  [V\otimes S_j:S_z]  
\end{pmatrix} = \begin{pmatrix}
  \chi_V(m_1)\chi_j(m_1) \\
  \vdots \\
  \chi_V(m_z)\chi_j(m_z)
\end{pmatrix}$$
where $1\le j \le z$ and $\chi_V$ and $\chi_j$ denote the characters for $V$ and $S_j$ respectively, and $m_1,\dots, m_z$ are representatives for the generalized conjugacy classes of $M$. Taking the $i$-th entry, we get: ($i$-th column of $X(M))$ $\cdot$ ($j$-th column of $M$) $=\chi_V(m_i)\chi_j(m_i)$,
but this just says that the $i$-th column of $X(M)$ is left eigenvector of $N(V)$ with eigenvalue $\chi_V(m_i)$. 
\end{proof}

For $1\le i\le m$, let $C_{i,1},\dots, C_{i,k_i}$ denote the conjugacy classes of $G_{e_i}$, and pick a complete set of representatives $g_{i,1},\dots, g_{i,k_i}$ for these classes. Let $V$ be a $kM$-module with character $\chi_V$. By \autoref{eigenvector lemma}, when $p\nmid \lvert G_{e_i}\rvert, 1\le i \le m$ each conjugacy class $C_{i,j}$ is canonically associated with a pair of left and right eigenvectors with eigenvalue $\chi_V(g_{i,j})$. We will refer to these as the $(i,j)$-th row of $(X(M))^T$ and the $(i,j)$-th column of $L^{-1}\overline{X}$, respectively. If $G_{e_i}$ is a maximal subgroup of $M$, extending notation from before we let $Z_V(G_{e_i})$ denote the elements $g \in G_{e_i}$ which act as scalars on $V$, and let $\omega_V(g)$ 
denote the corresponding scalar. We will write $k(n)$ for an asymptotic expression of $l(n)$.    

\begin{Theorem}\label{semisimple}Let $V$ be a $kM$ module with character $\chi_V$ and assume $p\nmid \lvert G_{e_i}\rvert,1\le i \le m$. 
Let $S_{i,j}$ denotes the sum over the $(i,j)$-th column of  $L^{-1}\overline{X}$, then  
\begin{enumerate} 

    \item the exact growth rate of $l(n)$ is 
     \begin{equation}\label{eq:exact}
l(n)=\sum_{i=1}^m \frac{1}{\lvert G_{e_i}\rvert }\sum_{j=1}^{k_i}\lvert C_{i,j}\rvert S_{i,j} \big( \chi(g_{i,j}) \big)^n         
     \end{equation}
    \item An asymptotic formula for $l(n)$ of $V$ is
\begin{equation}\label{eq:asym}
 l(n)\sim k(n):=\sum_{i=1}^m \frac{(\dim V)^n}{\lvert G_{e_i}\rvert }\sum_{j: g_{i,j} \in Z_V(G_{e_i})}\lvert C_{i,j}\rvert S_{i,j} \big( \omega_V(g_{i,j}) \big)^n.       
    \end{equation}
    \item We have $\lvert l(n)/k(n)-1\rvert \in \mathcal{O}\big((\lvert {\chi_{\mathrm{sec}}}\rvert / {\dim V})^n\big)$, where $\chi_{\mathrm{sec}}$ is any second largest character value (in terms of modulus) of $\chi_V$, and
    \item We have $\lvert l(n)-k(n)\rvert \in \mathcal{O}\big((\lvert \chi_{\mathrm{sec}}\rvert)^n\big)$.
\end{enumerate}
\end{Theorem}

\begin{proof}
For the eigenvalue $\chi(g_{i,j})$, choose the left eigenvector $w_{ij}^T$ to be the $(i,j)$-th row of $(X(M))^T$, \textcolor{black} {which as recalled above is the same as $X^T L$}, and choose the right eigenvector $v_{i,j}$ to be the $(i,j)$-th column of $L^{-1}\overline{X}$, normalized by a factor of $\lvert C_{i,j}\rvert /\lvert G_{e_i}\rvert$. From eigendecomposition we get
$$(N(V))^n=\sum_{i=1}^m w_{i,j}^Tv_{i,j} (\chi_V(g_{i,j}))^n.$$
We have $w_{i,j}^T v_{i,j}=1$ and if we sum over the column of $v_{i,j}w_{i,j}^T$ corresponding to the trivial $kM$-module we get $\lvert C_{i,j}\rvert S_{i,j}/\lvert G_{e_i}\rvert$. Statement (1) now follows as in \cite[Theorem 6]{lacabanne2023asymptotics}. Statements (2), (3) and (4) follow by noting that for $n$ large, the terms in the summation with $\lvert \chi_V(g_{i,j})\rvert = \dim V$ will dominate.     
\end{proof}

\begin{Remark}
If $M=G$ is a group, then if $p\nmid \lvert G\rvert$ \autoref{semisimple} gives $$l(n)=b(n) = \frac{1}{\vert G\vert} \sum_{t=1}^N \vert C_t\vert S_t \big(\chi(g_t)\big)^n$$ where for $1\le i \le N$, $t_i$ is a representative for the $i$-th conjugacy class $C_i$ of $G$, and $S_t=\overline{S_t}$ is the sum over entries in column of the character table corresponding to $g_t$. Thus \autoref{semisimple} generalizes the statement in \cite[Theorem 3]{he2024growthproblemsrepresentationsfinite}.    
\end{Remark}

\begin{Remark}
As we observed in the beginning of the section, if $kM$ is semisimple then $b(n)=l(n)$, so \autoref{semisimple} gives a refinement of the result of \autoref{Main} in this case. In particular, if no element other than $1$ act as identity on $M$, then the second statement of \autoref{semisimple} recovers the formula of $a(n)$ in \autoref{Main} in this case: this is because the submatrix of $L$ corresponding to the group of units $G$ is the identity matrix, so $S_{i,j}$ is exactly a column of the character table of $G$. Of course, if $kM$ is semisimple then we must have $p\nmid \lvert G_{e_i}\rvert$ for all $1\le i \le m$, so \autoref{semisimple} applies (see \textit{e.g.}\cite[Theorem 5.19]{steinberg2016representation}).    
\end{Remark}

\section{Examples}\label{examples}

\subsection{Full transformation monoids $T_n$, and symmetric inverse monoids $I_n$}

Let $T_n$ and $I_m$ denote the full transformation monoid and the symmetric inverse monoid on $m$ elements, respectively. In both cases the group of units is $G=S_m$. 

\begin{Proposition}\label{transformation monoids}
Suppose that either 
\begin{enumerate}
    \item $p =0$ and $M=T_m$, or
    \item $p \nmid m!$ and $M=I_m$.
\end{enumerate} Let $V$ be a $kM$-module on which no element other than 1 acts as identity. Then we have 
\begin{equation}\label{eqn:transformation_monoid}
b(n)\sim a(n)= \sum_{z=0}^{\floor * {m/2}} \frac{1}{(m-2z)!z!2^z}\cdot (\dim V)^n=k(n)\sim l(n). 
\end{equation}
\end{Proposition}

\begin{proof}
We know that $kI_m$ is semisimple under the assumptions of the theorem (see \cite[Corollary 9.4]{steinberg2016representation}). By \cite[Corollary A.3]{steinberg2016globaldimensiontransformationmonoid}, if $W$ is any $kS_m$-module which is not the sign representation, then $\Ind_G(W)$ is injective as a $kT_m$-module. Thus the assumptions of \autoref{Main} are satisfied for both $T_m$ and $I_m$, and the asymptotic formula for $S_m$ given in \cite[Example 2.3]{coulembier2023asymptotic} now implies the result for $a(n)$. The statement about $k(n)$ follows from \autoref{semisimple}.  
\end{proof}
\begin{Remark}
One might expect that for $T_m$ we can also weaken the requirement that $p=0$ to $p\nmid m!$. However, this is nontrivial to show, as monoid algebras need not have the same structure as we vary between characteristics that do not divide the order of any maximal subgroup. However, for $p$ sufficiently large the result of \autoref{transformation monoids} will hold.  
\end{Remark}

Since $kI_m$ is semisimple when $p\nmid m!$, the exact formula for $b(n)=l(n)$ can be calculated from \autoref{semisimple}. For example, if $V$ is the irreducible $\C I_3$-module with character values $0,1,2,0,3,1,0$, then using the character table and decomposition matrix for $\C I_3$ which can be found in \cite[Example 9.17]{steinberg2016representation}, in this case we obtain $$b(n)=2\cdot 3^{n-1} - 2^n+1.$$

\begin{figure}[ht]
    \centering
    \begin{minipage}{0.4\textwidth}
    \centering
    \includegraphics[width=\textwidth]{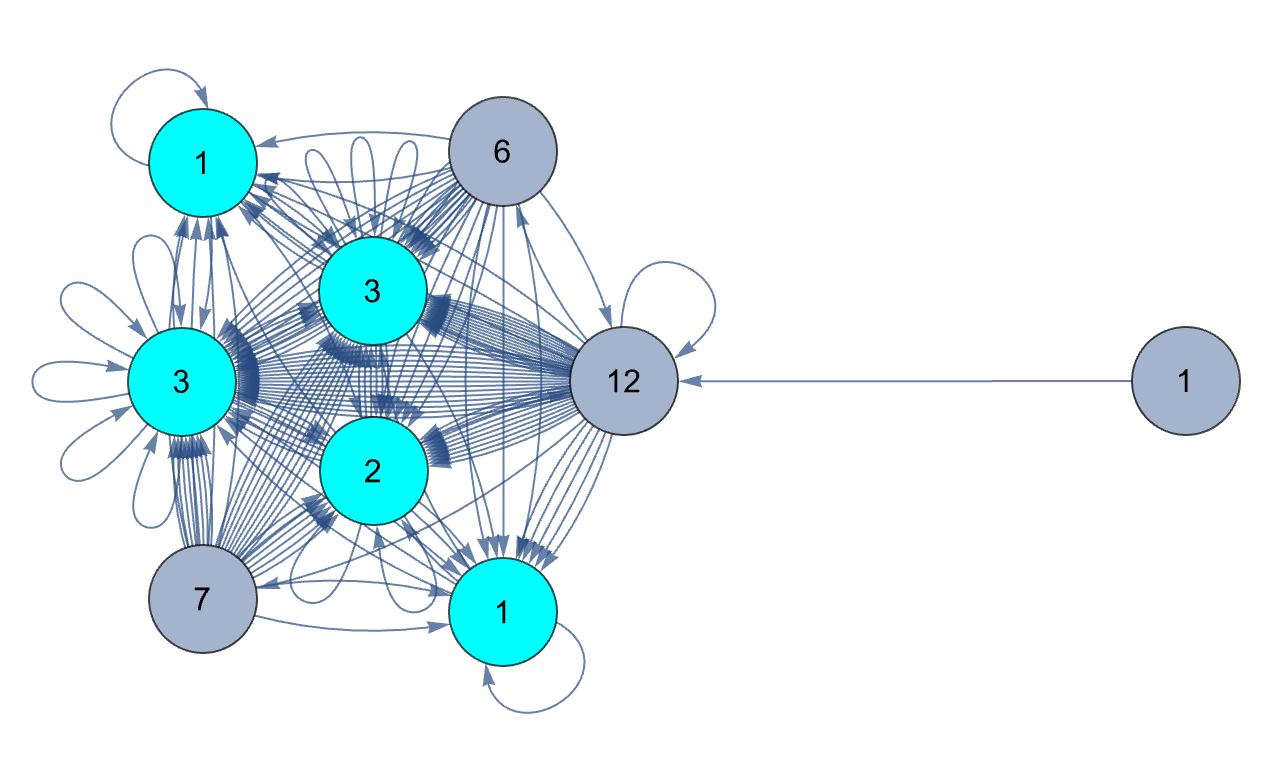} 
    \end{minipage}
    \begin{minipage}{0.3\textwidth}
     \centering
     
\[
\begin{psmallmatrix}
0 & 0 & 0 & 0 & 0 & 0 & 0 & 0 & 0 \\
1 & 1 & 0 & 0 & 0 & 0 & 0 & 1 & 1 \\
0 & 5 & 1 & 0 & 1 & 2 & 1 & 3 & 2 \\
0 & 4 & 0 & 1 & 1 & 1 & 2 & 2 & 4 \\
0 & 10 & 1 & 1 & 2 & 3 & 3 & 5 & 6 \\
0 & 15 & 2 & 1 & 3 & 5 & 4 & 8 & 8 \\
0 & 15 & 1 & 2 & 3 & 4 & 5 & 7 & 10 \\
0 & 1 & 0 & 0 & 0 & 0 & 0 & 0 & 0 \\
0 & 1 & 0 & 0 & 0 & 0 & 0 & 0 & 0
\end{psmallmatrix}
\]

    \end{minipage}
    \caption{The fusion graph and action matrix for the 12-dimensional projective $kT_4$-module in characteristic 0. The modules are labelled with their dimensions, and the projective cell $\Gamma^G$ is coloured in cyan.} 
    \label{fig:T4graph}
\end{figure}

\begin{figure}[H]
    \centering
    \includegraphics[scale=0.5]{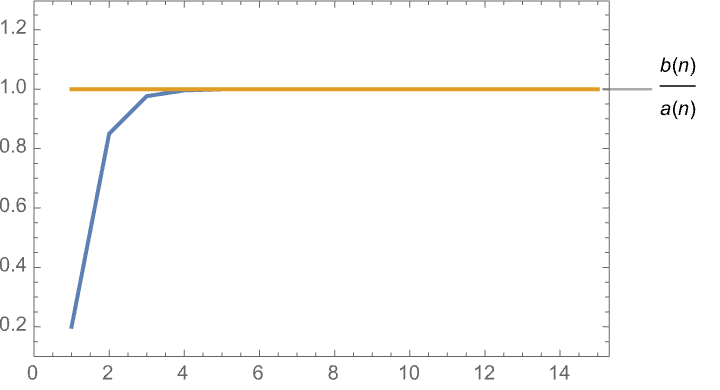}
    \caption{The ratio $b(n)/a(n)$ for the 12-dimensional projective $kT_4$-module in characteristic 0. Here $a(n)=5/12\cdot 12^n$ and $\lambda^{\textrm{sec}}=2$.}
    \label{fig:T4ratio}
\end{figure}

We plot below in \autoref{fig:T4graph} the fusion graph and action matrix for the 12-dimensional projective $kT_4$-module in characteristic 0, and in \autoref{fig:T4ratio} give the ratio $b(n)/a(n)$ in this case. In this case $\lambda^{\textrm{sec}}=2$, so $\lambda^{\textrm{sec}}/\dim V=1/6$, which explains the rapid convergence. \textcolor{black}{We provide python code in \cite{code} to generate the multiplication table of $T_m$ and then set up the regular representation of $kT_m$ in Magma.} As usual, we approximated the characteristic zero behavior by doing the computation in a large enough finite field. 

\subsection{Matrix monoids $\M(2,q)$}

Let $\M(2,q)=\M(2,\mathbb{F}_q)$ denote the monoid of $2\times 2$  matrices with entries in $\mathbb{F}_q$, and let $G=\GL(2,q)$ be its group of units. We know from \cite{kovacs1992semigroup} that when $p \nmid \lvert G\rvert = (q^2-1)(q^2-q)$, $k\M(2,q)$ is semisimple. If $q$ is odd, let $\epsilon$ denote the unique element of order 2 in the center $Z(G)$ of $G$.

\begin{Proposition}[Matrix monoids $\M(2,q)$ in nonmodular characteristic]
Suppose $p\nmid \lvert G\rvert$ and $V$ is a $k\M(2,q)$-module. We have \begin{equation}\label{nonmodular matrix monoid}
a(n)=\begin{cases}
\frac{q^2+(-1)^n}{(q+1)(q^2-q)} \cdot (\dim V)^n & \text{if $\epsilon$ acts as scalar on $\Res_G(V)$}, \vspace{5pt} \\
\frac{q^2}{(q+1)(q^2-q)} \cdot (\dim V)^n & \text{else}.
\end{cases} 
\end{equation} 
If $q$ is even, then we are always in the second case.
\end{Proposition}

\begin{proof}
Since $\M(2,q)$ is semisimple, by \autoref{Main} (or \autoref{semisimple}) in this case $a(n)$ is determined by the character table of $G$, which can be found in \textit{e.g.} \cite[\S 5.2]{fulton2013representation}. The sum over the dimensions of irreducible $kG$-modules is $q^2(q-1)$, and when $q$ is odd the sum over the column of the character table corresponding to $\epsilon$ is $q-1$. The sum over all other columns in the center vanishes, yielding our formula.     
\end{proof}

Next we consider $\M(2,q)$ in defining characteristic.

\begin{Proposition}[Matrix monoids $\M(2,q)$ in defining characteristic]
Let $q=p^r$ be a prime power where $p$ is also the characteristic of $k$.  If $V$ is a $k\M(2,q)$-module on which no element other than 1 acts as identity, then \begin{equation}\label{eqn: conjecture}
a(n)=\begin{cases}
\frac{(p+1)^r}{2^r(q+1)(q-1)}\Big(1+ \frac{1}{q}(-1)^n\Big)\cdot (\dim V)^n &  \text{if $\epsilon$ acts as scalar on $\Res_G(V)$}\\
\frac{(p+1)^r}{2^r(q+1)(q-1)} \cdot (\dim V)^n & \text{else.}
\end{cases}    
\end{equation}
If $q$ is even, then we are always in the second case.
\end{Proposition}

\begin{proof}
Let $G=\textup{GL}(2,q)$ be the group of units of $kM(2,q)$. By \cite[Theorem 4 and Theorem 12]{KOUWENHOVEN1993369}, most of the projective indecomposable $kG$-modules are injective when considered as $k\M(2,q)$-modules. So condition ii) in \autoref{Main} is satisfied, and it suffices to compute the asymptotic formula for a faithful $kG$-module. Thus, the statement of the proposition follows from \cite[Theorem 13]{he2024growthproblemsrepresentationsfinite}. 
\end{proof}

\begin{figure}[H]
    \centering
\includegraphics[scale=0.5]{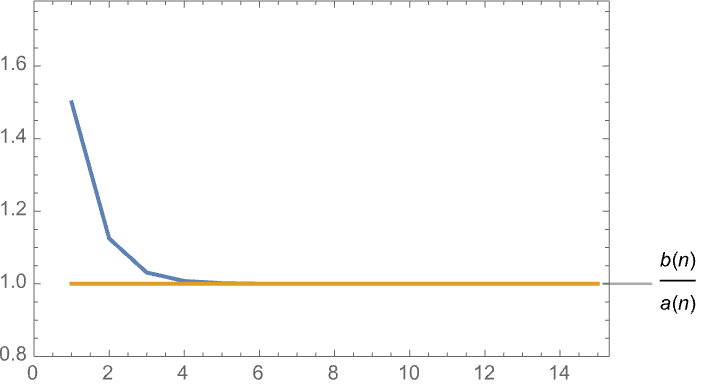}
    \caption{The ratio $b(n)/a(n)$ for the four-dimensional projective $k\M(2,3)$-module in characteristic 3. Here $a(n)=(1/4+1/12(-1)^n)\cdot 4^n$ and $\lambda^{\text{sec}} = 1.$}
    \label{fig:M23ratio}
\end{figure}

We give in \autoref{fig:km(2,3)} below the fusion graph and action matrix for the four-dimensional projective $k\M(2,3)$-module in defining characteristic. We also plot the ratio $b(n)/a(n)$ in this case in \autoref{fig:M23ratio}. In this case a second largest eigenvalue is $\lambda^{\sec}=1$.

\begin{figure}[ht]
    \centering
    \begin{minipage}{0.40\textwidth}
    \centering
\includegraphics[width=\textwidth]{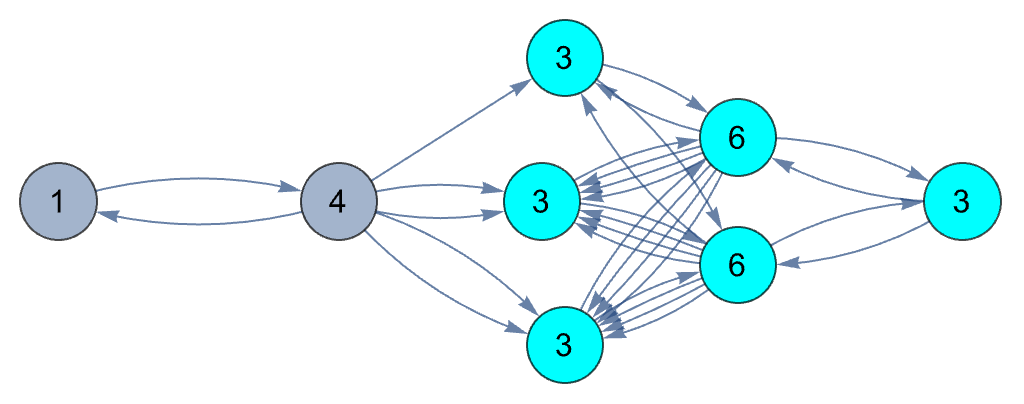} 
    \end{minipage}
    \begin{minipage}{0.4\textwidth}
     \centering
     
\[
\begin{psmallmatrix}
0 & 1 & 0 & 0 & 0 & 0 & 0 & 0 \\
1 & 0 & 0 & 0 & 0 & 0 & 0 & 0 \\
0 & 2 & 0 & 0 & 0 & 3 & 3 & 0 \\
0 & 1 & 0 & 0 & 0 & 1 & 1 & 0 \\
0 & 2 & 0 & 0 & 0 & 3 & 3 & 0 \\
0 & 0 & 1 & 1 & 1 & 0 & 0 & 1 \\
0 & 0 & 1 & 1 & 1 & 0 & 0 & 1 \\
0 & 0 & 0 & 0 & 0 & 1 & 1 & 0
\end{psmallmatrix}
\]

    \end{minipage}
    \caption{The fusion graph and action matrix for the four-dimensional projective $k\M(2,3)$-module in characteristic 3. The modules are labelled with their dimensions, and the projective cell $\Gamma^G$ is coloured in cyan.} 
    \label{fig:km(2,3)}
\end{figure}

\subsection{An example where \autoref{Main} does not apply}\label{SS:Counterexample}

Condition ii) in the statement of \autoref{Main} is somewhat restrictive: it is already not satisfied when the monoid is $N=\{1,x,0\}$ where $x^2=0$, and of the 35 nonisomorphic monoids of order 4, 9 do not satisfy the condition in characteristic 0. Generalizing the example of $N$, if $M$ is any monoid with group of units $G$, the direct product $M\times N$ (where $N$ may be replaced with any nontrivial nilpotent monoid) will not satisfy condition ii), as any $kG$-module gets a non-split self-extension from the nilpotent element $x$ and so cannot be injective. (These counterexamples were communicated to us by Walter Mazorchuk.)

 However, experimentation with monoids of small order (see \autoref{conj section}) gives evidence that $V$, when condition i) is satisfied, will still reach $\Gamma^P$ even when ii) is not satisfied. As an example, let $M=T_3\times N$ where $T_3$ is the full transformation monoid on 3 elements and $N=\{1,x,0\}$ as above. There is a unique four-dimensional $kM$-module $V$ satisfying condition i), in characteristic 31. In this case we get $a(n)=2/3\cdot 4^n$ as predicted by \autoref{eqn:asym} despite condition ii) not being fulfilled. The fusion graph and action matrix for $V$ are given below, showing that the projective cell is still reached.

\begin{figure}[ht]
    \centering
    \begin{minipage}{0.40\textwidth}
    \centering
    \includegraphics[width=\textwidth]{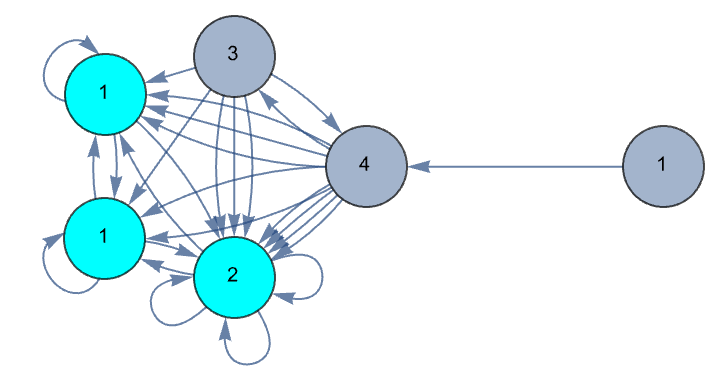} 
    \end{minipage}
    \begin{minipage}{0.4\textwidth}
     \centering
     
\[
\begin{pmatrix}
0 & 0 & 0 & 0 & 0 & 0 \\
1 & 0 & 0 & 0 & 0 & 1 \\
0 & 3 & 1 & 1 & 1 & 1 \\
0 & 2 & 1 & 1 & 1 & 1 \\
0 & 4 & 1 & 1 & 3 & 3 \\
0 & 1 & 0 & 0 & 0 & 0
\end{pmatrix}
\]

    \end{minipage}
    \caption{The fusion graph and action matrix for the  four-dimensional projective $k(T_3\times N)$-module $V$, in characteristic 31. The modules are labelled with their dimensions, and the projective cell $\Gamma^G$ is coloured in cyan.} 
    \label{fig:T3Ngraph}
\end{figure}






\begin{thebibliography}{CEOT24}

\bibitem[BK72]{bryant1972tensor}
R.M. Bryant and L.G. Kov{\'a}cs.
\newblock Tensor products of representations of finite groups.
\newblock {\em Bull. London Math. Soc.}, 4:133--135, 1972.
\newblock \href {https://doi.org/10.1112/blms/4.2.133}
  {\path{doi:10.1112/blms/4.2.133}}.

\bibitem[CEO24]{coulembier2023asymptotic}
K.~Coulembier, P.~Etingof, and V.~Ostrik.
\newblock Asymptotic properties of tensor powers in symmetric tensor
  categories.
\newblock {\em Pure Appl. Math. Q.}, 20(3):1141--1179, 2024.
\newblock URL: \url{https://arxiv.org/abs/2301.09804}, \href
  {https://doi.org/10.4310/pamq.2024.v20.n3.a4}
  {\path{doi:10.4310/pamq.2024.v20.n3.a4}}.

\bibitem[CEOT25]{coulembier2024fractalbehaviortensorpowers}
K.~Coulembier, P.~Etingof, V.~Ostrik, and D.~Tubbenhauer.
\newblock Fractal behavior of tensor powers of the two dimensional space in
prime characteristic.
\newblock {\em Contemp. Math.}, 
Volume: 829; 2025; 291 pp, 
American Mathematical Society, [Providence], RI, [2025], ©2025.
\newblock \url{https://arxiv.org/abs/2405.16786}, \href
{https://doi.org/10.1090/conm/829}
{\path{doi:10.1090/conm/829}}.

\bibitem[COT24]{Coulembier_2023}
K.~Coulembier, V.~Ostrik, and D.~Tubbenhauer.
\newblock Growth rates of the number of indecomposable summands in tensor
  powers.
\newblock {\em Algebr. Represent. Theory}, 27(2):1033--1062, 2024.
\newblock URL: \url{https://arxiv.org/abs/2301.00885}, \href
  {https://doi.org/10.1007/s10468-023-10245-7}
  {\path{doi:10.1007/s10468-023-10245-7}}.

\bibitem[FH91]{fulton2013representation}
W.~Fulton and J.~Harris.
\newblock {\em Representation theory}, volume 129 of {\em Graduate Texts in
  Mathematics}.
\newblock Springer-Verlag, New York, 1991.
\newblock A first course, Readings in Mathematics.
\newblock \href {https://doi.org/10.1007/978-1-4612-0979-9}
  {\path{doi:10.1007/978-1-4612-0979-9}}.

\bibitem[He25]{he2024growthproblemsrepresentationsfinite}
D.~He.
\newblock Growth problems for representations of finite groups.
\newblock 2025. To appear in Ark. Mat.
\newblock URL: \url{https://arxiv.org/abs/2408.04196}.

\bibitem[HT25]{code}
D.~He and D.~Tubbenhauer.
\newblock Growth problems for representations of finite monoids, 2025.
\newblock URL: \url{https://github.com/dtubbenhauer/monoidgrowth}.

\bibitem[KST24]{KhSiTu-monoidal-cryptography}
M.~Khovanov, M.~Sitaraman, and D.~Tubbenhauer.
\newblock Monoidal categories, representation gap and cryptography.
\newblock {\em Trans. Amer. Math. Soc. Ser. B}, 11:329--395, 2024.
\newblock URL: \url{https://arxiv.org/abs/2201.01805}, \href
  {https://doi.org/10.1090/btran/151} {\path{doi:10.1090/btran/151}}.

\bibitem[Kou93]{KOUWENHOVEN1993369}
F.M. Kouwenhoven.
\newblock Indecomposable representations of {$M(2,\mathbf F_q)$} over {$\mathbf
  F_q$}.
\newblock {\em J. Algebra}, 155(2):369--396, 1993.
\newblock \href {https://doi.org/10.1006/jabr.1993.1050}
  {\path{doi:10.1006/jabr.1993.1050}}.

\bibitem[Kov92]{kovacs1992semigroup}
L.G. Kov\'acs.
\newblock Semigroup algebras of the full matrix semigroup over a finite field.
\newblock {\em Proc. Amer. Math. Soc.}, 116(4):911--919, 1992.
\newblock \href {https://doi.org/10.2307/2159467} {\path{doi:10.2307/2159467}}.

\bibitem[LTV23]{lacabanne2023asymptotics}
A.~Lacabanne, D.~Tubbenhauer, and P.~Vaz.
\newblock Asymptotics in finite monoidal categories.
\newblock {\em Proc. Amer. Math. Soc. Ser. B}, 10:398--412, 2023.
\newblock URL: \url{https://arxiv.org/abs/2307.03044}, \href
  {https://doi.org/10.1090/bproc/198} {\path{doi:10.1090/bproc/198}}.

\bibitem[LTV24]{lacabanne2024asymptotics}
A.~Lacabanne, D.~Tubbenhauer, and P.~Vaz.
\newblock Asymptotics in infinite monoidal categories.
\newblock 2024.
\newblock To appear in High. Struct.
\newblock \url{https://arxiv.org/abs/2404.09513}.

\bibitem[LPRS24]{lachowska2024tensorpowersvectorrepresentation}
A.~Lachowska, O.~Postnova, N.~Reshetikhin, and D.~Solovyev.
\newblock Tensor powers of vector representation of $u_q(\mathfrak{sl}_2)$ at
  even roots of unity.
\newblock 2024.
\newblock URL: \url{https://arxiv.org/abs/2404.03933}.

\bibitem[Lar24]{larsen2024boundsmathrmsl2indecomposablestensorpowers}
M.J. Larsen.
\newblock Bounds for $\mathrm{SL}_2$-indecomposables in tensor powers of the
  natural representation in characteristic $2$.
\newblock 2024.
\newblock URL: \url{https://arxiv.org/abs/2405.16015}.

\bibitem[Ste16b]{steinberg2016representation}
B.~Steinberg.
\newblock {\em Representation theory of finite monoids}.
\newblock Universitext. Springer, Cham, 2016.
\newblock \href {https://doi.org/10.1007/978-3-319-43932-7}
  {\path{doi:10.1007/978-3-319-43932-7}}.

\bibitem[Ste16a]{steinberg2016globaldimensiontransformationmonoid}
B.~Steinberg.
\newblock The global dimension of the full transformation monoid (with an
  appendix by {V}. {M}azorchuk and {B}. {S}teinberg).
\newblock {\em Algebr. Represent. Theory}, 19(3):731--747, 2016.
\newblock URL: \url{https://arxiv.org/abs/1502.00959}, \href
  {https://doi.org/10.1007/s10468-016-9597-4}
  {\path{doi:10.1007/s10468-016-9597-4}}.

\bibitem[Ste23]{steinberg2023modularrepresentationtheorymonoids}
B.~Steinberg.
\newblock The modular representation theory of monoids.
\newblock 2023.
\newblock URL: \url{https://arxiv.org/abs/2305.08251}.

\end{thebibliography}
\end{document}